# Hidden Convexity-Based Distributed Operation of Integrated Electricity-Gas Systems

Rong-Peng Liu, *Member, IEEE*, Yue Song, *Member, IEEE*, Junhong Liu, *Student Member, IEEE*, Xiaozhe Wang, *Senior Member, IEEE*, Jinpeng Guo, *Member, IEEE*, and Yunhe Hou, *Senior Member, IEEE*

*Abstract*—This letter proposes a hidden convexity-based method to address distributed optimal energy flow (OEF) problems for transmission-level integrated electricity-gas systems. First, we develop a node-wise decoupling method to decompose an OEF problem into multiple OEF subproblems. Then, we propose a hidden convexity-based method to equivalently reformulate nonconvex OEF subproblems as semi-definite programs. This method differs from any approximation and convexification methods that may incur infeasible solutions. Since all OEF subproblems are originally convex or equivalently convexified, we adopt an ADMM to solve the hidden convexity-based distributed OEF problem with convergence analysis. Test results validate the effectiveness of the proposed method, especially in handling a large number of agents.

## NOMENCLATURE

### A. Sets

| | |
|---|---|
| $\mathcal{D}_n^P/\mathcal{D}_n^G$ | Set of power/gas loads connecting to power/gas node $n$. |
| $\mathcal{L}_n^{P+}/\mathcal{L}_n^{G+}/\mathcal{C}_n^{G+}$ | Set of power transmission lines/gas passive pipelines/gas compressors connecting to power/gas/gas node $n$ with a presupposed inflow direction to node $n$. ($\mathcal{L}_n^{P-}/\mathcal{L}_n^{G-}/\mathcal{C}_n^{G-}$ is with a presupposed outflow direction from node $n$.) |
| $\mathcal{N}_n^P/\mathcal{N}_n^G$ | Set of (neighboring) power/gas nodes that directly connect to power/gas node $n$. |
| $\mathcal{P}/\mathcal{G}$ | Set of power/gas nodes. |
| $\mathcal{U}_n^P/\mathcal{U}_{n\text{-}m}^G/\mathcal{W}_n^G$ | Set of coal- & gas-fired units/gas-fired units/gas wells connecting to power node $n$/power node $n$ and gas node $m$/gas node $n$. |

### B. Parameters

| | |
|---|---|
| $C_g(\cdot)/C_w$ | Cost function of coal-fired unit $g$/cost of gas well $w$. |
| $G_n^{\min}/G_n^{\max}$ | Pressure square limit of gas node $n$. |
| $G_w^{\min}/G_w^{\max}$ | Output limit of gas well $w$. |
| $P_d/G_d$ | Power/gas load $d$. |
| $P_g^{\min}/P_g^{\max}$ | Output limit of unit $g$. |
| $P_{mn}/G_{mn}/C_{mn}$ | Capacity of power transmission line $mn$/gas passive pipeline $mn$/gas compressor $mn$. |
| $x_{mn}/W_{mn}/\alpha_{mn}$ | Reactance of power transmission line $mn$/Weymouth constant of gas passive pipeline $mn$/ratio of gas compressor $mn$. |
| $\theta_n^{\min}/\theta_n^{\max}$ | Phase angle limit of power node $n$. |
| $\chi_g$ | Conversion ratio of gas-fired unit $g$. |

### C. Variables

| | |
|---|---|
| $p_g/g_w$ | Output of unit $g$/gas well $w$. |
| $p_{mn}/g_{mn}/c_{mn}$ | Power/gas/gas flow in power transmission line $mn$/gas passive pipeline $mn$/gas compressor $mn$ with a predefined flow direction from node $m$ to $n$. |
| $\theta_n/\pi_n$ | Phase angle/pressure square of power node $n$/gas node $n$. |
| $\lambda$ | Lagrangian multiplier (vector). |

## I. INTRODUCTION

AS is reported by the U.S. Energy Information Administration [1], in the U.S., natural gas outdistanced other energy sources for electricity generation in 2020 and was predicted to possibly remain a bellwether in 2050, unveiling its competitiveness. In order to raise the transmission efficiency of natural gas to gas-fired power generators and reduce the transmission cost, integrated electricity-gas systems (IEGSs) have been constructed massively [2]. As a result, the optimal operation of IEGSs, especially the IEGS optimal energy flow (OEF) problems, becomes a research hotspot.

IEGS OEF problems study (short-term) optimal IEGS operation strategies for some specific objective(s), e.g., operation costs, under IEGS operation constraints [2]-[8]. Reference [2] focuses on centralized OEF problems. In view of the multi-agent nature [3]-[5] and increased emphasis on utility data privacy [6]-[9], privacy-preserving distributed operation of IEGSs becomes attractive, leading to *distributed OEF problems*. However, nonconvex gas flow constraints, e.g., bidirectional Weymouth equations (mixed-integer nonconvex equations), incur significant challenges for solving distributed OEF problems. Previous works [3]-[8] propose various methods for handling Weymouth equations in distributed OEF problems, including piecewise linear approximation [3], second-order cone relaxation [4], extended convex hull relaxation [5], and convex-concave procedure (CCP) [6]-[8]. As is analysed in [6] and [7], approximation- and relaxation-based methods [3]-[5] cannot ensure the feasibility of solutions. For CCP, a local heuristic method, its solution quality relies on initial points. Even though, the convergence of CCP can only be guaranteed to the critical points of original problems (but not necessarily local optimal solutions) [10]. In addition, previous works [3]-

This work was supported in part by the Fonds de recherche du Québec-Nature et technologies (FRQ-NT) under Grant 334636, in part by the Fonds de Recherche du Quebec-Nature et technologies under Grant FRQ-NT PR-298827, FRQ-NT 2023-NOVA-314338, and in part by the National Natural Science Foundation of China under Grant 52177188. *(Corresponding author: Junhong Liu.)*

R. Liu and X. Wang are with the Department of Electrical and Computer Engineering, McGill University, Montreal, QC H3A 0E9, Canada (e-mail: rpliu@eee.hku.hk/rongpeng.liu@mail.mcgill.ca; xiaozhe.wang2@mcgill.ca).

Y. Song is with the State Key Laboratory of Intelligent Autonomous Systems and Frontiers Science Center for Intelligent Autonomous Systems, Tongji University, Shanghai, China (e-mail: ysong@tongji.edu.cn).

J. Liu and Y. Hou are with the Department of Electrical and Electronic Engineering, The University of Hong Kong, Hong Kong SAR, China (e-mail: jhliu@eee.hku.hk; yhhou@eee.hku.hk).

J. Guo is with the Department of Energy and Electrical Engineering, Ho-Hai University, Nanjing, Jiangsu, China. (email: jinpeng.guo@hhu.edu.cn).

[8] consider only a limited number of agents ($\leq 4$). The performances of their proposed methods (for handling Weymouth equations) on multiple agents in distributed OEF problems are under investigation.

By exploring the structure features of bidirectional Weymouth equations in distributed OEF problems, we find that they are hidden convex [11]. Based on this property, this letter addresses the challenge in solving multi-agent distributed OEF problems caused by Weymouth equations. First, we develop a node-wise decoupling method for decomposing an OEF problem into multiple OEF subproblems. Then, we propose a hidden convexity-based method to eliminate the nonconvexity in the (nonconvex) OEF subproblems that contain Weymouth equations. Accordingly, these (nonconvex) subproblems are *equivalently* reformulated as semi-definite programs (SDPs). Since all subproblems are originally convex or equivalently convexified, we adopt an ADMM [12] to solve hidden convexity-based distributed OEF problems with convergence analysis. Test results validate the effectiveness of the proposed method, especially in handling *a large number of agents*.

Note that the hidden convexity-based method proposed in [13] (for addressing the DistFlow model that has only continuous variables) cannot be directly adopted to cope with Weymouth equations due to the existence of binary variables, as strong duality does not hold for general optimization problems with discrete variables. To the best of the authors' knowledge, this is *the first work* to address distributed OEF problems with (original) bidirectional Weymouth equations, completely differing from approximation- and relaxation-based methods [3]-[5] that may derive infeasible solutions. In addition, unlike CCP [6]-[8], the hidden convexity-based method does not rely on initial points and is ensured to converge to the *global optimums of the nonconvex subproblems*.

The remainder of this letter is organized as follows. Section II presents the distributed OEF formulation. Section III proposes the hidden convexity-based solution method. Case studies are conducted in Section IV. Section V draws conclusions.

## II. Decomposed OEF Formulation

In order to enable the distributed operation of IEGSs, this section develops a node-wise decoupling method by decomposing a (centralized) OEF problem into an equivalent distributed OEF problem that consists of multiple OEF subproblems. Before proceeding, we clarify the scope of this work.

1) This work focuses on OEF problems for *transmission-level* IEGSs and thus adopts Weymouth equations and simplified gas compressor constraints for modelling gas flow in gas passive and active pipelines, respectively [2]-[8].

2) This work focuses on single-period OEF problems and thus does not consider i) energy storage facilities and line pack of gas pipelines and ii) different time scales of power and gas systems. These important topics remain our future work.

*Notations*: Scalars are denoted by regular lowercase or uppercase letters. Vectors and matrices are denoted by bold lowercase and uppercase letters, respectively. All vectors are column vectors. Superscript $^T$ denotes a transpose manipulation.

First, we introduce the formulation of a centralized OEF problem, denoted as $\mathbf{P_o}$.

$\mathbf{P_o}$: 
$$\min \sum_{n \in \mathcal{P}} \sum_{g \in \mathcal{U}_n^P \setminus \mathcal{U}_{n-m}^G} C_g(p_{g(n)}) + \sum_{n \in \mathcal{G}} \sum_{w \in \mathcal{W}_n^G} C_w \cdot g_{w(n)} \quad (1)$$

s.t. $P_g^{\min} \leq p_{g(n)} \leq P_g^{\max}, \quad \forall g \in \mathcal{U}_n^P,$ (2a)

$\theta_n^{\min} \leq \theta_n \leq \theta_n^{\max},$ (2b)

$-P_{mn} \leq p_{mn(n)} \leq P_{mn}, \quad \forall mn \in \mathcal{L}_n^{P+},$ (2c)

$x_{mn} \cdot p_{mn(n)} = \theta_{m(n)} - \theta_{n(n)}, \quad \forall mn \in \mathcal{L}_n^{P+},$ (2d)

$\sum_{g \in \mathcal{U}_n^P} p_{g(n)} + \sum_{mn \in \mathcal{L}_n^{P+}} p_{mn(n)} - \sum_{no \in \mathcal{L}_n^{P-}} p_{no(n)} = \sum_{d \in \mathcal{D}_n^P} P_{d(n)},$ (2e)

$\forall n \in \mathcal{P}$ in (2a)-(2e), (2f)

$G_w^{\min} \leq g_{w(n)} \leq G_w^{\max}, \quad \forall w \in \mathcal{W}_n^G,$ (3a)

$G_n^{\min} \leq \pi_n \leq G_n^{\max},$ (3b)

$0 \leq g_{mn(n)} \leq G_{mn}, \quad \forall mn \in \mathcal{L}_n^{G+},$ (3c)

$0 \leq c_{mn(n)} \leq C_{mn}, \quad \forall mn \in \mathcal{C}_n^{G+},$ (3d)

$(g_{mn(n)})^2 = W_{mn}^2 \cdot (\pi_{m(n)} - \pi_{n(n)}), \quad \forall mn \in \mathcal{L}_n^{G+},$ (3e)

$\pi_{n(n)} \leq \alpha_{mn} \cdot \pi_{m(n)}, \quad \forall mn \in \mathcal{C}_n^{G+},$ (3f)

$\sum_{w \in \mathcal{W}_n^G} g_{w(n)} + \sum_{mn \in \mathcal{L}_n^{G+}} g_{mn(n)} - \sum_{no \in \mathcal{L}_n^{G-}} g_{no(n)} + \sum_{in \in \mathcal{C}_n^{G+}} c_{mn(n)} - \sum_{nj \in \mathcal{C}_n^{G-}} c_{nj(n)} = \sum_{d \in \mathcal{D}_n^G} G_{d(n)} + \sum_{g \in \mathcal{U}_{m-n}^P} \chi_g \cdot p_{g(n)},$ (3g)

$\forall n \in \mathcal{G}$ in (3a)-(3g). (3h)

Subscript $x(n)$ denotes IEGS component $x$ that connects to power/gas node $n$. Objective (1) minimizes IEGS operation costs, i.e., power generator and gas well costs. Function $C_g(\cdot)$ is a quadratic function with a non-negative quadratic coefficient. Constraints (2) are power subsystem constraints, including generator output, phase angle, power transmission line capacity, DC power flow [14], and power balance constraints. Constraints (3) are gas subsystem constraints, including gas well output, nodal pressure square, gas passive pipeline capacity, gas compressor capacity, unidirectional Weymouth equation [6], gas compressor [3]-[7], and gas balance constraints.

Then, we proceed to the distributed OEF formulation. By implementing the node-wise decoupling method on $\mathbf{P_o}$, we derive multiple OEF subproblems. Specifically, for each power node $n, n \in \mathcal{P}$, we introduce pseudo variables $\{p_{g(n)}^c, p_{mn(n)}^c, p_{no(n)}^c, \theta_{r(n)}^c\}, g \in \mathcal{U}_n^P, mn \in \mathcal{L}_n^{P+}, no \in \mathcal{L}_n^{P-}, r \in \mathcal{N}_n^P \cup \{n\}$. Based on the node-wise decoupling method, the OEF subproblem for *individual* power node $n$ is

$$\min_{\mathbf{x}_n^P, \mathbf{y}_n^P} \sum_{g \in \mathcal{U}_n^P \setminus \mathcal{U}_{n-m}^G} C_g(p_{g(n)}) \quad (4a)$$

s.t. constraints (2a)-(2c), (4b)

$\begin{cases} x_{mn} \cdot p_{mn(n)}^c = \theta_{m(n)}^c - \theta_{n(n)}^c, & \forall mn \in \mathcal{L}_n^{P+} \\ x_{no} \cdot p_{no(n)}^c = \theta_{n(n)}^c - \theta_{o(n)}^c, & \forall no \in \mathcal{L}_n^{P-}, \end{cases}$ (4c)

$\sum_{g \in \mathcal{U}_n^P} p_{g(n)}^c + \sum_{mn \in \mathcal{L}_n^{P+}} p_{mn(n)}^c - \sum_{no \in \mathcal{L}_n^{P-}} p_{no(n)}^c = \sum_{d \in \mathcal{D}_n^P} P_{d(n)}.$ (4d)

$\mathbf{x}_n^P = col(p_{g(n)}, p_{mn(n)}, \theta_n), g \in \mathcal{U}_n^P, mn \in \mathcal{L}_n^{P+}$, and $\mathbf{y}_n^P = col(p_{g(n)}^c, p_{mn(n)}^c, p_{no(n)}^c, \theta_{r(n)}^c), g \in \mathcal{U}_n^P, mn \in \mathcal{L}_n^{P+}, no \in \mathcal{L}_n^{P-}, r \in \mathcal{N}_n^P \cup \{n\}$. Function $col(\cdot)$ reshapes scalars and/or vectors as one vector, e.g., $col(\mathbf{a}, b, \mathbf{c}) = [\mathbf{a}^T \, b \, \mathbf{c}^T]^T$.

For each gas node $n, n \in \mathcal{G}$, we introduce pseudo variables $\{g_{w(n)}^c, g_{mn(n)}^c, g_{no(n)}^c, c_{in(n)}^c, c_{nj(n)}^c, \pi_{r(n)}^c, p_{g(n)}'^c\}, w \in \mathcal{W}_n^G, mn \in \mathcal{L}_n^{G+}, no \in \mathcal{L}_n^{G-}, in \in \mathcal{C}_n^{G+}, nj \in \mathcal{C}_n^{G-}, r \in \mathcal{N}_n^G \cup \{n\}, g \in \mathcal{U}_{m-n}^G$.

Based on the node-wise decoupling method, the OEF subproblem for *individual* gas node $n$ is

$$\min_{\mathbf{x}_n^G, \mathbf{y}_n^G} \sum_{w \in \mathcal{W}_n^G} C_w \cdot g_{w(n)} \quad (5a)$$

s.t. constraints (3a)-(3d), (5b)

$$\begin{cases} g_{mn(n)}^c = \sqrt{W_{mn}^2 \cdot (\pi_{m(n)}^c - \pi_{n(n)}^c)}, & \forall mn \in \mathcal{L}_n^{G+} \\ g_{no(n)}^c = \sqrt{W_{no}^2 \cdot (\pi_{n(n)}^c - \pi_{o(n)}^c)}, & \forall no \in \mathcal{L}_n^{G-}, \end{cases} \quad (5c)$$

$$\begin{cases} \pi_{n(n)}^c \le \alpha_{mn} \cdot \pi_{m(n)}^c, & \forall mn \in \mathcal{C}_n^{G+} \\ \pi_{o(n)}^c \le \alpha_{no} \cdot \pi_{n(n)}^c, & \forall no \in \mathcal{C}_n^{G-}, \end{cases} \quad (5d)$$

$$\sum_{w \in \mathcal{W}_n^G} g_{w(n)}^c + \sum_{mn \in \mathcal{L}_n^{G+}} g_{mn(n)}^c - \sum_{no \in \mathcal{L}_n^{G-}} g_{no(n)}^c + \sum_{in \in \mathcal{C}_n^{G+}} c_{in(n)}^c - \sum_{nj \in \mathcal{C}_n^{G-}} c_{nj(n)}^c = \sum_{d \in \mathcal{D}_n^G} G_{d(n)} + \sum_{g \in \mathcal{U}_{m-n}^G} \chi_g \cdot p'^c_{g(n)}. \quad (5e)$$

$\mathbf{x}_n^G = col(g_{w(n)}, g_{mn(n)}, c_{in(n)}, \pi_n)$, $w \in \mathcal{W}_n^G$, $mn \in \mathcal{L}_n^{G+}$, $in \in \mathcal{C}_n^{G+}$, and $\mathbf{y}_n^G = col(g_{w(n)}^c, g_{mn(n)}^c, g_{no(n)}^c, c_{in(n)}^c, c_{nj(n)}^c, \pi_{r(n)}^c, p'^c_{g(n)})$, $w \in \mathcal{W}_n^G$, $mn \in \mathcal{L}_n^{G+}$, $no \in \mathcal{L}_n^{G-}$, $in \in \mathcal{C}_n^{G+}$, $nj \in \mathcal{C}_n^{G-}$, $r \in \mathcal{N}_n^G \cup \{n\}$, $g \in \mathcal{U}_{m-n}^G$.

To ensure the equivalence between the (centralized) OEF problem $\mathbf{P_o}$ (before decoupling) and all OEF subproblems (4)-(5), all $n \in \mathcal{P}$ in (4) and all $n \in \mathcal{G}$ in (5) (after decoupling), pseudo variables should reach a consensus with their counterparts, yielding the following coupling constraints.

$$p_{g(n)}^c = p_{g(n)}, \quad \forall g \in \mathcal{U}_n^P, n \in \mathcal{P}, \quad (6a)$$

$$\theta_{n(i)}^c = \theta_n, \quad \forall i \in \mathcal{N}_n^P \cup \{n\}, n \in \mathcal{P}, \quad (6b)$$

$$p_{mn(m)}^c = p_{mn(n)}, p_{mn(n)}^c = p_{mn(n)}, \quad \forall mn \in \mathcal{L}^{P+}, n \in \mathcal{P}, \quad (6c)$$

$$g_{w(n)}^c = g_{w(n)}, \quad \forall w \in \mathcal{W}_n^G, n \in \mathcal{G}, \quad (7a)$$

$$\pi_{n(m)}^c = \pi_n, \quad \forall m \in \mathcal{N}_n^G \cup \{n\}, n \in \mathcal{G}, \quad (7b)$$

$$g_{mn(m)}^c = g_{mn(n)}, g_{mn(n)}^c = g_{mn(n)}, \quad \forall mn \in \mathcal{L}_n^{G+}, n \in \mathcal{G}, \quad (7c)$$

$$c_{mn(m)}^c = c_{mn(n)}, c_{mn(n)}^c = c_{mn(n)}, \quad \forall mn \in \mathcal{C}_n^{G+}, n \in \mathcal{G}, \quad (7d)$$

$$p'^c_{g(n)} = p_{g(m)}, \quad \forall g \in \mathcal{U}_{m-n}^G, n \in \mathcal{G}. \quad (7e)$$

***Remark 1:*** In practice, an IEGS may incorporate multiple agents (rather than only two power and gas subsystem agents). For example, based on geographical partition, even a connected power/gas subsystem may belong to several entities [5], [7]. The node-wise decoupling method, rendering each node in an IEGS an individual agent, is one of the most flexible partitioning methods and, meanwhile, preserves the privacy of all entities to the greatest extent. In addition, the widely employed optical-fiber communication network is considered for IEGS data transmission, which is resilient against adverse interceptions and further ensures the privacy of individual agents [15]. Note that this decoupling method is still applicable even if an entity consists of multiple nodes without any modification. So far, by means of pseudo variables, the *distributed OEF problem* (4)-(7), all $n \in \mathcal{P}$ in (4) and all $n \in \mathcal{G}$ in (5), is equivalent to centralized OEF problem $\mathbf{P_o}$. The only challenge for solving the distributed OEF problem (by distributed methods) is the nonconvex Weymouth equation (5c).

## III. HIDDEN CONVEXITY-BASED DISTRIBUTED OPERATION

In this section, we first propose a hidden convexity-based method for eliminating nonconvexity in the OEF subproblems that contain nonconvex Weymouth equations. Then, we present a distributed solution method for solving the hidden convexity-based distributed OEF problem.

### A. Hidden Convexity-Based Nonconvexity Elimination

For gas passive pipelines $mn(n)$ and $no(n)$, $mn \in \mathcal{L}_n^{G+}$, $no \in \mathcal{L}_n^{G-}$, we introduce pseudo variables $\mathbf{y}_{mn(n)}^C = col(g'^c_{mn(n)}, \pi'^c_{m(n)}, \pi'^c_{n_m(n)})$, $\mathbf{y}_{no(n)}^C = col(g'^c_{no(n)}, \pi'^c_{n_o(n)}, \pi'^c_{o(n)})$, $mn \in \mathcal{L}_n^{G+}$, $no \in \mathcal{L}_n^{G-}$. Accordingly, we rewrite nonconvex (5c) as

$$\begin{cases} g_{mn(n)}^c \le \sqrt{W_{mn}^2 \cdot (\pi_{m(n)}^c - \pi_{n(n)}^c)}, & \forall mn \in \mathcal{L}_n^{G+} \\ g_{no(n)}^c \le \sqrt{W_{no}^2 \cdot (\pi_{n(n)}^c - \pi_{o(n)}^c)}, & \forall no \in \mathcal{L}_n^{G-}, \end{cases} \quad (8a)$$

$$\begin{cases} g'^c_{mn(n)} \ge \sqrt{W_{mn}^2 \cdot (\pi'^c_{m(n)} - \pi'^c_{n_m(n)})}, & \forall mn \in \mathcal{L}_n^{G+} \\ g'^c_{no(n)} \ge \sqrt{W_{no}^2 \cdot (\pi'^c_{n_o(n)} - \pi'^c_{o(n)})}, & \forall no \in \mathcal{L}_n^{G-}. \end{cases} \quad (8b)$$

To ensure the equivalence between (5c), all $n \in \mathcal{G}$ in (5c), and (8), all $n \in \mathcal{G}$ in (8), pseudo variables are enforced by

$$g'^c_{mn(m)} = g_{mn(n)}, g'^c_{mn(n)} = g_{mn(n)}, \quad \forall mn \in \mathcal{L}_n^{G+}, n \in \mathcal{G}, \quad (9a)$$

$$\pi'^c_{n(m)} = \pi_n, \pi'^c_{n_m(m)} = \pi_n, \quad \forall m \in \mathcal{N}_n^G, n \in \mathcal{G}. \quad (9b)$$

For an *entire* IEGS, the compact form of a distributed OEF problem is

$$\min_{\mathbf{x}, \mathbf{y}} \sum_{n \in \mathcal{P}} f_n(\mathbf{x}_n^P) + \sum_{n \in \mathcal{G}} g_n(\mathbf{x}_n^G) \quad (10a)$$

s.t. $\mathbf{Ax} = \mathbf{By}$, (10b)

$$\mathbf{x}_n^P \in \mathcal{X}_n^P, \mathbf{y}_n^P \in \mathcal{Y}_n^P, \quad \forall n \in \mathcal{P}, \quad (10c)$$

$$\mathbf{x}_n^G \in \mathcal{X}_n^G, \mathbf{y}_n^G \in \mathcal{Y}_n^G, \quad \forall n \in \mathcal{G}, \quad (10d)$$

$$\mathbf{y}_{mn(n)}^C \in \mathcal{Y}_{mn(n)}^C, \mathbf{y}_{no(n)}^C \in \mathcal{Y}_{no(n)}^C,$$
$$\forall mn \in \mathcal{L}_n^{G+}, no \in \mathcal{L}_n^{G-}, n \in \mathcal{G}. \quad (10e)$$

$\mathbf{x} = col(\mathbf{x}_{n'}^P, \mathbf{x}_n^G)$ and $\mathbf{y} = col(\mathbf{y}_{n'}^P, \mathbf{y}_n^G, \mathbf{y}_{mn(n)}^C, \mathbf{y}_{no(n)}^C)$, $mn \in \mathcal{L}_n^{G+}$, $no \in \mathcal{L}_n^{G-}$, $n' \in \mathcal{P}, n \in \mathcal{G}$. Functions $f_n(\cdot)$ and $g_n(\cdot)$ in (10a) denote objective functions (4a) and (5a), respectively. Constraint (10b) is the compact form of coupling constraints (6), (7), and (9), all $n \in \mathcal{G}$ in (9). $\mathbf{A}$ and $\mathbf{B}$ are constant matrices. In (10c)-(10e), sets $\mathcal{X}_n^P = \{\mathbf{x}_n^P | (4b)\}$, $\mathcal{X}_n^G = \{\mathbf{x}_n^G | (5b)\}$, $\mathcal{Y}_n^P = \{\mathbf{y}_n^P | (4c), (4d)\}$, $\mathcal{Y}_n^G = \{\mathbf{y}_n^G | (5d), (5e), (8a)\}$, $\mathcal{Y}_{mn(n)}^C = \{\mathbf{y}_{mn(n)}^C | (8b)\}$, and $\mathcal{Y}_{no(n)}^C = \{\mathbf{y}_{no(n)}^C | (8b)\}$. Now, all the sets in (10c)-(10e) consist of only local variables. Considering *each set as an individual agent*, problem (10) can be iteratively solved by the ADMM [12] *in parallel*, i.e.,

$$\mathbf{x}_n^{P\,k+1}/\mathbf{x}_n^{G\,k+1} = \underset{\mathbf{x}_n^P \in \mathcal{X}_n^P / \mathbf{x}_n^G \in \mathcal{X}_n^G}{\arg\min} \mathcal{L}_d(\mathbf{x}, \mathbf{y}^k, \boldsymbol{\lambda}^k), \quad \forall n \in \mathcal{P} / \forall n \in \mathcal{G}, \quad (11a)$$

$$\mathbf{y}_n^{P\,k+1}/\mathbf{y}_n^{G\,k+1}/\mathbf{y}_{mn(n)}^{C\,k+1}/\mathbf{y}_{no(n)}^{C\,k+1} = \underset{\substack{\mathbf{y}_n^P \in \mathcal{Y}_n^P / \mathbf{y}_n^G \in \mathcal{Y}_n^G / \mathbf{y}_{mn(n)}^C \in \\ \mathcal{Y}_{mn(n)}^C / \mathbf{y}_{no(n)}^C \in \mathcal{Y}_{no(n)}^C}}{\arg\min} \mathcal{L}_d(\mathbf{x}^{k+1}, \mathbf{y}, \boldsymbol{\lambda}^k),$$

$$\forall n \in \mathcal{P} / \forall n \in \mathcal{G} / \forall mn \in \mathcal{L}_n^{G+}, n \in \mathcal{G} / \forall no \in \mathcal{L}_n^{G-}, n \in \mathcal{G}, \quad (11b)$$

$$\boldsymbol{\lambda}_n^{P\,k+1}/\boldsymbol{\lambda}_n^{G\,k+1}/\boldsymbol{\lambda}_{mn(n)}^{C\,k+1}/\boldsymbol{\lambda}_{no(n)}^{C\,k+1} = \boldsymbol{\lambda}_n^{P\,k}/\boldsymbol{\lambda}_n^{G\,k}/\boldsymbol{\lambda}_{mn(n)}^{C\,k}/\boldsymbol{\lambda}_{no(n)}^{C\,k} + d \circ (\mathbf{A}_n \mathbf{x}_n^{P\,k+1} /$$
$$\mathbf{A}_n \mathbf{x}_n^{G\,k+1} - \mathbf{B}_n \mathbf{y}_n^{P\,k+1} / \mathbf{B}_n \mathbf{y}_n^{G\,k+1} / \mathbf{B}_{mn(n)} \mathbf{y}_{mn(n)}^{C\,k+1} / \mathbf{B}_{no(n)} \mathbf{y}_{no(n)}^{C\,k+1}),$$

$$\forall n \in \mathcal{P} / \forall n \in \mathcal{G} / \forall mn \in \mathcal{L}_n^{G+}, n \in \mathcal{G} / \forall no \in \mathcal{L}_n^{G-}, n \in \mathcal{G}, \quad (11c)$$

where $\mathcal{L}_d$ is the augmented Lagrangian function,

$$\mathcal{L}_d(\mathbf{x}, \mathbf{y}, \boldsymbol{\lambda}) = h(\mathbf{x}) + \boldsymbol{\lambda}^T (\mathbf{Ax} - \mathbf{By}) + (1/2) \cdot d \cdot \|\mathbf{Ax} - \mathbf{By}\|_2^2. \quad (12)$$

The superscript $k/k+1$ in (11) denotes the $k^{th}/(k+1)^{th}$ ADMM iteration. Scalar d in (11c) is a constant. Matrices $\mathbf{A}_n$, $\mathbf{B}_n$, $\mathbf{B}_{mn(n)}$, and $\mathbf{B}_{no(n)}$ are constant matrices. Symbol $\circ$ refers to the pointwise production of a scalar and a vector/matrix. Function $h(\cdot)$ in (12) denotes objective function (10a). The Lagrangian multiplier $\boldsymbol{\lambda} = col(\boldsymbol{\lambda}_{n'}^{P\,k+1}, \boldsymbol{\lambda}_n^{G\,k+1}, \boldsymbol{\lambda}_{mn(n)}^{C\,k+1}, \boldsymbol{\lambda}_{no(n)}^{C\,k+1})$, $mn \in \mathcal{L}_n^{G+}$, $no \in \mathcal{L}_n^{G-}$, $n' \in \mathcal{P}, n \in \mathcal{G}$.

Unfortunately, feasible regions $\mathcal{Y}_{mn(n)}^{C}$ and $\mathcal{Y}_{no(n)}^{C}$, $mn \in \mathcal{L}_n^{G+}$, $no \in \mathcal{L}_n^{G-}$, $n \in \mathcal{G}$, in (11b) are nonconvex due to (8b), challenging the implementation of the ADMM. In fact, the (nonconvex) subproblems in (11b) that contain $\mathcal{Y}_{mn(n)}^{C}$ and $\mathcal{Y}_{no(n)}^{C}$ are hidden-convex. We take gas passive pipeline $mn(n)$ as an example (the same story for $no(n)$), and the compact form (of its corresponding nonconvex subproblem in (11b)) is

$$\min_{\mathbf{y}_{mn(n)}^{C}} (1/2) \cdot \mathbf{y}_{mn(n)}^{C\,T} \mathbf{A}_{mn(n)}^{0} \mathbf{y}_{mn(n)}^{C} + \mathbf{b}_{mn(n)}^{0\,T} \mathbf{y}_{mn(n)}^{C} + c_{mn(n)}^{0} \quad (13a)$$

$$\text{s.t.} (1/2) \cdot \mathbf{y}_{mn(n)}^{C\,T} \mathbf{A}_{mn(n)}^{1} \mathbf{y}_{mn(n)}^{C} + \mathbf{b}_{mn(n)}^{1\,T} \mathbf{y}_{mn(n)}^{C} \leq 0. \quad (v_{mn(n)}) \quad (13b)$$

$\mathbf{A}_{mn(n)}^{0}$ and $\mathbf{A}_{mn(n)}^{1}$ are constant matrices. $\mathbf{b}_{mn(n)}^{0}$ and $\mathbf{b}_{mn(n)}^{1}$ are constant vectors. $c_{mn(n)}^{0}$ is a constant. $v_{mn(n)}$ is the dual variable for constraint (13b). Based on the duality theory [11], the dual problem for (13) is

$$\max_{v_{mn(n)}} -(1/2) \cdot [(\mathbf{b}_{mn(n)}^{0} + v_{mn(n)} \circ \mathbf{b}_{mn(n)}^{1})^T (\mathbf{A}_{mn(n)}^{0} + v_{mn(n)} \circ \mathbf{A}_{mn(n)}^{1})^{-1} \cdot$$
$$(\mathbf{b}_{mn(n)}^{0} + v_{mn(n)}^{*} \circ \mathbf{b}_{mn(n)}^{1}) + c_{mn(n)}^{0}] \quad (14a)$$

$$\text{s.t. } \mathbf{A}_{mn(n)}^{0} + v_{mn(n)} \circ \mathbf{A}_{mn(n)}^{1} \succeq 0, \quad (14b)$$
$$v_{mn(n)} \geq 0. \quad (14c)$$

By implementing the Schur complement [11], we re-write (14) as the following SDP.

$$\max_{v_{mn(n)}, \gamma_{mn(n)}} \gamma_{mn(n)} \quad (15a)$$

$$\text{s.t. } v_{mn(n)} \geq 0, \quad (15b)$$

$$\begin{bmatrix} \mathbf{A}_{mn(n)}^{0} + v_{mn(n)} \circ \mathbf{A}_{mn(n)}^{1} & \mathbf{b}_{mn(n)}^{0} + v_{mn(n)} \circ \mathbf{b}_{mn(n)}^{1} \\ (\mathbf{b}_{mn(n)}^{0} + v_{mn(n)} \circ \mathbf{b}_{mn(n)}^{1})^T & 2 \cdot (c_{mn(n)}^{0} - \gamma_{mn(n)}) \end{bmatrix} \succeq 0, \quad (15c)$$

where $\gamma_{mn(n)}$ is an auxiliary variable.

***Remark 2***: Strong duality between (13) and (14) holds if the feasible region $\mathcal{Y}_{mn(n)}^{C}$ is non-empty, as the optimization problem (13) is a quadratically constrained quadratic program with one and only one (quadratic) constraint (i.e, QCQP-1) [11]. For this case, $\mathbf{y}_{mn(n)}^{C} = col(0, 0, 0)$ is always a feasible solution. Thus, the duality gap between (13) and (14) (as well as (15)) is zero, i.e., the strong duality holds. The same story for $\mathcal{Y}_{no(n)}^{C}$. In general, for any gas passive pipeline $mn(n)/no(n)$, $mn \in \mathcal{L}_n^{G+}/no \in \mathcal{L}_n^{G-}$, $n \in \mathcal{G}$, nonconvexity in (11b) can be eliminated by replacing each nonconvex subproblem (in (11b)) with (15). Thus, we can derive *global optimums* of nonconvex subproblems in (11b) by solving their SDP counterparts (instead of solving nonconvex subproblems in (11b)). So far, all OEF subproblems in (11a) and (11b) are either originally convex or equivalently convexified (as SDPs), enabling the implementation of the ADMM. Note that the proposed hidden convexity-based (nonconvexity elimination) method only applies to hidden convex optimization problems, such as nonconvex subproblems in (11b) (i.e., QCQP-1). In other words, this method cannot be directly applied to general nonconvex OEF problems. For example, OEF subproblem (5) is not hidden convex, and this method cannot equivalently reformulate (5) as a convex optimization problem.

### B. Bidirectional Weymouth Equation

In order to focus on illustrating the hidden convexity-based method, the above analyze only considers (relatively simpler) unidirectional Weymouth equation (5c), a nonconvex quadratic equation. This subsection shows how the hidden convexity-based method is employed to address the (original) bidirectional Weymouth equation, a mixed-integer nonconvex quadratic equation. For gas passive pipeline $mn(n)$, the (original) bidirectional Weymouth equation is

$$\text{sgn}(\pi_{m(n)}^{c}, \pi_{n(n)}^{c}) \cdot g_{mn(n)}^{c} = \sqrt{\text{sgn}(\pi_{m(n)}^{c}, \pi_{n(n)}^{c}) \cdot W_{mn}^{2} \cdot (\pi_{m(n)}^{c} - \pi_{n(n)}^{c})}, \quad (16a)$$

$$\text{sgn}(\pi_{m(n)}^{c}, \pi_{n(n)}^{c}) = 1/-1, \quad \pi_{m(n)}^{c} \geq \pi_{n(n)}^{c} / \pi_{m(n)}^{c} < \pi_{n(n)}^{c}. \quad (16b)$$

By introducing a constant $G_{\max} = \max\{G_{mn}\}$, $mn \in \mathcal{L}_n^{G+}, n \in \mathcal{G}$, an auxiliary variable $u_{mn(n)}$, model (16) is equivalent to

$$(g_{mn(n)}^{c})^{2} = W_{mn}^{2} \cdot (\pi_{m(n)}^{c} - \pi_{n(n)}^{c}) u_{mn(n)}, \quad (17a)$$

$$(u_{mn(n)} + 1)(u_{mn(n)} - 1) = 0, \quad (17b)$$

$$G_{\max} \cdot (u_{mn(n)} - 1) \leq g_{mn(n)}^{c} \leq G_{\max} \cdot (u_{mn(n)} + 1). \quad (17c)$$

Constraints (17) are nonconvex due to *quadratic equations* (17a) and (17b). Similar to (5c), they are hidden-convex constraints. Their nonconvexity can be eliminated by the proposed method in Section III.A by the following procedures: i) splitting (17a) and (17b) into two in-equality constraints, respectively, by means of pseudo variables (similar to (8)); ii) building coupling relations between pseudo variables and the counterparts of their original variables (similar to (9)); iii) constructing subproblems, each for one inequality constraint in Step i) (similar to the subproblems in (11b)); iv) reformulating the nonconvex subproblems constructed in Step iii) as SDPs based on the (strong) duality theorem. The same story for gas passive pipeline $no(n)$. Note that gas passive pipeline transmission constraint (3c) should be replaced by (18) when considering bidirectional Weymouth equation (16).

$$-G_{mn} \leq g_{mn(n)} \leq G_{mn}, \quad \forall mn \in \mathcal{L}_n^{G+}, \quad (18)$$

### C. Solution Method

We adopt the ADMM [11] to solve the hidden convexity-based distributed OEF problem, and details are presented in Algorithm 1.

---

**Algorithm 1**: ADMM for solving the hidden convexity-based distributed OEF problem

---

1: **Initialize** variables $\mathbf{x}^0$, $\mathbf{y}^0$, and $\boldsymbol{\lambda}^0$, parameter d (d > 0), and convergence thresholds $\varepsilon_{\text{pri}}$ and $\varepsilon_{\text{dual}}$, and set $k = 0$;
2: **while** $k = 0$ or either $(20a) \leq \varepsilon_{\text{pri}}$ or $(20b) \leq \varepsilon_{\text{dual}}$ is not satisfied, **do**
3:   Solve (11a) in parallel. Update $\mathbf{x}^{k+1}$;
4:   Solve (11b) in parallel by replacing each nonconvex subproblem in (11b) with (15). Derive $v_{mn(n)}^{*}, v_{no(n)}^{*}$. Update $\mathbf{y}^{k+1}$, where

$$\mathbf{y}_{mn(n)}^{C\,k+1} = -(\mathbf{A}_{mn(n)}^{0} + v_{mn(n)}^{*} \circ \mathbf{A}_{mn(n)}^{1})^{-1} (\mathbf{b}_{mn(n)}^{0} + v_{mn(n)}^{*} \circ \mathbf{b}_{mn(n)}^{1}), \quad (19a)$$

$$\mathbf{y}_{no(n)}^{C\,k+1} = -(\mathbf{A}_{no(n)}^{0} + v_{no(n)}^{*} \circ \mathbf{A}_{no(n)}^{1})^{-1} (\mathbf{b}_{no(n)}^{0} + v_{no(n)}^{*} \circ \mathbf{b}_{no(n)}^{1}); \quad (19b)$$

---

5:     Update $\boldsymbol{\lambda}_{n'}^{P\,k+1}$, $\boldsymbol{\lambda}_n^{G\,k+1}$, $\boldsymbol{\lambda}_{mn(n)}^{C\,k+1}$, and $\boldsymbol{\lambda}_{no(n)}^{C\,k+1}$ for $mn \in \mathcal{L}_n^{G+}$, $no \in \mathcal{L}_n^{G-}$, $n' \in \mathcal{P}$, $n \in \mathcal{G}$, in parallel based on (11c);

6:     Set $k = k+1$;

7:     Calculate primal and dual residuals

$$\|\mathbf{Ax}^k - \mathbf{By}^k\|_2, \quad (20a)$$

$$d \cdot \|\mathbf{A}^T\mathbf{B}(\mathbf{y}^k - \mathbf{y}^{k-1})\|_2, \quad (20b)$$

where $\|\cdot\|_2$ denotes the $l_2$-norm;

8: **end while**

---

***Remark 3***: Equations (19a) and (19b) are derived based on KKT conditions (i.e., the stationary condition [11] of primal problem (13)). At each iteration, each agent solves his/her own subproblems in (11a), (11b), and (15) in turn and then updates his/her own Lagrangian multiplier in (11c) with partial operation data shared among neighboring agents (without any centralized coordinators). Algorithm 1 is iteratively implemented as indicated above until satisfying the convergence conditions in Step 2. Since the implementation of Algorithm 1 does not rely on any centralized coordinators, it empowers fully distributed operation of hidden convexity-based distributed OEF problems. Next, we analyze its convergence.

***Proposition 1***: Algorithm 1 (for solving the hidden convexity-based distributed OEF problem) converges to the KKT point of the (centralized) OEF problem $\mathbf{P_o}$ if either $(\mathbf{Ax}^k - \mathbf{By}^k) \to 0$ as $k \to \infty$ or $(\boldsymbol{\lambda}^{k+1} - \boldsymbol{\lambda}^k) \to 0$ as $k \to \infty$.

We first prove the following Lemma:

***Lemma 1***: Conditions (21) and (22) are equivalent for Algorithm 1.

$$(\mathbf{Ax}^k - \mathbf{By}^k) \to 0 \text{ as } k \to \infty. \quad (21)$$

$$(\boldsymbol{\lambda}^{k+1} - \boldsymbol{\lambda}^k) \to 0 \text{ as } k \to \infty. \quad (22)$$

***Proof***: For ease of reading, we present the compact form of updating rule (11) as follows:

$$\mathbf{x}^{k+1} = \arg \min_{\mathbf{x} \in \mathcal{X}} \mathcal{L}_d(\mathbf{x}, \mathbf{y}^k, \boldsymbol{\lambda}^k), \quad (23a)$$

$$\mathbf{y}^{k+1} = \arg \min_{\mathbf{y} \in \mathcal{Y}} \mathcal{L}_d(\mathbf{x}^{k+1}, \mathbf{y}, \boldsymbol{\lambda}^k), \quad (23b)$$

$$\boldsymbol{\lambda}^{k+1} = \boldsymbol{\lambda}^k + d \circ (\mathbf{Ax}^{k+1} - \mathbf{By}^{k+1}). \quad (23c)$$

The definitions of vectors $\mathbf{x}$, $\mathbf{y}$, and $\boldsymbol{\lambda}$ are given in the second paragraph of Section III.A. Set $\mathcal{X}$ is the union of all $\mathcal{X}_{n'}^P$ and $\mathcal{X}_n^G$, $n' \in \mathcal{P}$, $n \in \mathcal{G}$, and set $\mathcal{Y}$ is the union of all $\mathcal{Y}_{n'}^P$, $\mathcal{Y}_n^G$, $\mathcal{Y}_{mn(n)}^C$, and $\mathcal{Y}_{no(n)}^C$, $mn \in \mathcal{L}_n^{G+}$, $no \in \mathcal{L}_n^{G-}$, $n' \in \mathcal{P}$, $n \in \mathcal{G}$.

If condition (21) holds, according to (23c), we have

$$(\boldsymbol{\lambda}^{k+1} - \boldsymbol{\lambda}^k) = d \circ (\mathbf{Ax}^{k+1} - \mathbf{By}^{k+1}) \to 0 \text{ as } k \to \infty,$$

i.e., condition (22) holds. Similarly, we can derive (21) if (22) holds. This completes the proof. ∎

Then, we prove Proposition 1 (inspired by [16]).

***Proof***: For ease of reading, we present the compact form of the (centralized) OEF problem $\mathbf{P_o}$ as follows.

$$\mathbf{P_c}: \min_{\mathbf{x}} \mathbf{f}^T\mathbf{x} \quad (24a)$$

$$\text{s.t. } \mathbf{Cx} = \mathbf{c}, \quad (24b)$$

$$\mathbf{x}^T\mathbf{D}_i\mathbf{x} = \mathbf{d}_i^T\mathbf{x}, \quad i = 1, \cdots, k, \quad (24c)$$

$$\mathbf{Ex} \le \mathbf{e}, \quad (24d)$$

$$\mathbf{x} \in \mathbb{R}^n. \quad (24e)$$

Based on the definition of $\mathbf{x}$, it is also the same variables in the centralized OEF problem $\mathbf{P_c}$. The Lagrangian function for $\mathbf{P_c}$ is

$$\mathcal{L}_P(\mathbf{x}, \boldsymbol{\gamma}, \boldsymbol{\eta}, \boldsymbol{\mu}) = \mathbf{f}^T\mathbf{x} + \boldsymbol{\gamma}^T(\mathbf{Cx} - \mathbf{c}) +$$

$$\sum_{i=1}^{k} \eta_i (\mathbf{x}^T\mathbf{D}_i\mathbf{x} - \mathbf{d}_i^T\mathbf{x}) + \boldsymbol{\mu}^T(\mathbf{Ex} - \mathbf{e}), \quad (25)$$

where $\boldsymbol{\gamma}$, $\boldsymbol{\eta}$ ($\boldsymbol{\eta} = col(\eta_i)$, $i = 1, \cdots, k$), and $\boldsymbol{\mu}$ are Lagrangian multipliers.

Assuming that i) $\mathbf{x}^*$ is a (local) optimal solution for $\mathbf{P_c}$, and ii) $\mathbf{P_c}$ satisfies some mild conditions [11], according to the KKT conditions [11], we have $\boldsymbol{\gamma}^*$, $\boldsymbol{\eta}^*$, and $\boldsymbol{\mu}^*$ (for centralized OEF problem $\mathbf{P_c}$), such that

$$\mathbf{f}^T + (\boldsymbol{\gamma}^*)^T\mathbf{C} + \sum_{i=1}^{k} \eta_i^* (2(\mathbf{x}^*)^T\mathbf{D}_i - \mathbf{d}_i^T) + (\boldsymbol{\mu}^*)^T\mathbf{E} = \mathbf{0}, \quad (26a)$$

$$\mathbf{Cx}^* = \mathbf{c}, \quad (26b)$$

$$(\mathbf{x}^*)^T\mathbf{D}_i\mathbf{x}^* = \mathbf{d}_i^T\mathbf{x}^*, \quad i = 1, \cdots, k, \quad (26c)$$

$$\mathbf{Ex}^* \le \mathbf{e}, \quad (26d)$$

$$\boldsymbol{\mu}^* \ge \mathbf{0}, \quad (26e)$$

$$(\boldsymbol{\mu}^*)^T(\mathbf{Ex}^* - \mathbf{e}) = 0. \quad (26f)$$

According to Lemma 1 and the assumptions in Proposition 1, we have i) $(\mathbf{Ax}^k - \mathbf{By}^k) \to 0$ and ii) $(\boldsymbol{\lambda}^k - \boldsymbol{\lambda}^{k-1}) \to 0$. Inspired by [16], based on (12) and the updating rule (23), the $(k+1)^{th}$ iteration in Algorithm 1 satisfies

$$\mathbf{f}^T + (\boldsymbol{\lambda}^k)^T\mathbf{A} + d \circ (\mathbf{Ax}^{k+1} - \mathbf{By}^k)^T\mathbf{A} + \boldsymbol{\varsigma}_1^T\mathbf{E}_1 = \mathbf{0}, \quad (27a)$$

$$-(\boldsymbol{\lambda}^k)^T\mathbf{B} - d \circ (\mathbf{Ax}^{k+1} - \mathbf{By}^{k+1})^T\mathbf{B} + \boldsymbol{\varrho}^T\mathbf{C}' +$$

$$\sum_{i=1}^{k} \sigma_i (2(\mathbf{y}^{k+1})^T\mathbf{D}_i' - (\mathbf{d}_i')^T) + \boldsymbol{\varsigma}_2^T\mathbf{E}_2' = \mathbf{0}. \quad (27b)$$

$\boldsymbol{\varsigma}_1$, $\boldsymbol{\varsigma}_2$, $\boldsymbol{\varrho}$, and $\boldsymbol{\sigma}$ ($\boldsymbol{\sigma} = col(\sigma_i)$, $i = 1, \cdots, k$) are Lagrangian multipliers. $\mathbf{E}_1$, $\mathbf{E}_2'$, $\mathbf{C}'$, and $\mathbf{D}_i'$, $i = 1, \cdots, k$, are constant matrices, and $\mathbf{d}_i'$, $i = 1, \cdots, k$, are constant vectors. we have $\mathbf{E} = COL[\mathbf{E}_1, \mathbf{E}_2'\mathbf{A}]$, $\mathbf{C} = \mathbf{C}'\mathbf{A}$, $\mathbf{D}_i = \mathbf{A}^T\mathbf{D}_i'\mathbf{A}$, and $\mathbf{d}_i = \mathbf{d}_i'\mathbf{A}$. Function $COL(\cdot)$ reshapes matrices as one matrix. In addition, it is easy to know that the matrix $\mathbf{B}$ is an identity matrix. Based on (23c) equations (27) are simplified as

$$\mathbf{f}^T + (\boldsymbol{\lambda}^{k+1})^T\mathbf{A} + \boldsymbol{\varsigma}_1^T\mathbf{E}_1 = \mathbf{0}, \quad (28a)$$

$$(\boldsymbol{\lambda}^{k+1})^T = \boldsymbol{\varrho}^T\mathbf{C}' + \sum_{i=1}^{k} \sigma_i (2(\mathbf{y}^{k+1})^T\mathbf{D}_i' - (\mathbf{d}_i')^T) + \boldsymbol{\varsigma}_2^T\mathbf{E}_2'. \quad (28b)$$

By integrating (28b) into (28a), we have

$$\mathbf{f}^T + \boldsymbol{\varrho}^T\mathbf{C} + \sum_{i=1}^{k} \sigma_i (2(\mathbf{y}^{k+1})^T\mathbf{D}_i'\mathbf{A} - \mathbf{d}_i^T) + [col(\boldsymbol{\varsigma}_1, \boldsymbol{\varsigma}_2)]^T\mathbf{E} = 0. \quad (29)$$

Since $(\mathbf{Ax}^k - \mathbf{By}^k) \to 0$, we can replace the $\mathbf{By}^{k+1}$ in (29) with $\mathbf{Ax}^{k+1}$. In addition, by setting $\boldsymbol{\varrho} = \boldsymbol{\gamma}^*$, $\boldsymbol{\sigma} = \boldsymbol{\eta}^*$, and $col(\boldsymbol{\varsigma}_1, \boldsymbol{\varsigma}_2) = \boldsymbol{\mu}^*$, equation (29) is equivalent to the KKT condition (26a). This result indicates that $\mathbf{x}^{k+1}$ (derived by Algorithm 1) satisfies KKT condition (26a). Meanwhile, this $\mathbf{x}^{k+1}$ also satisfies KKT conditions (26b)-(26f), as i) $\mathbf{x}^{k+1}$ and $\mathbf{y}^{k+1}$ satisfy constraints (24b)-(24e) and ii) $(\mathbf{Ax}^k - \mathbf{By}^k) \to 0$. Thus, KKT conditions (26) are all satisfied, indicating that Algorithm 1 converges to the KKT point of the (centralized) OEF problem $\mathbf{P_c}$. This completes the proof. ∎

## IV. CASE STUDY

To validate the effectiveness of the proposed method, this section conducts tests on a small-scale 6-node-power-7-node-gas system (IEGS-6-7) and a medium-scale 118-node-power-20-node gas system (IEGS-118-20). Detailed IEGS topology and parameters as well as the parameters in Algorithm 1 are presented in [17]. Codes are written in MATLAB and execut-ed on an i5-6500 PC. Mosek, Gurobi, IPOPT, and SCIP are solvers for SDPs, mixed-integer convex programs, continuous convex and nonconvex programs, and mixed-integer nonconvex programs, respectively. (See [17], [18] for more details.)

First, we compare the solutions of distributed OEF problems with different methods for handling Weymouth equations, i.e., piecewise linear approximation (PWL) [3], second-order cone relaxation (SOC) [4], CCP [7], and the proposed hidden convexity-based method (HCM), on the small-scale IEGS-6-7. Test results are shown in Table I. A vanilla ADMM [12] with Gurobi and IPOPT is adopted to solve PWL-, SOC- and CCP-based distributed OEF problems. Algorithm 1 with Mosek and IPOPT is used to solve HCM-based distributed OEF problems. |Opt. gap| = |(Obj. – Obj. of NCW)/(Obj. of NCW)|. Function |·| denotes the absolute value. The solution for the centralized OEF problem (solved by IPOPT and SCIP) based on nonconvex Weymouth equations (NCW) is used as a reference to validate solution feasibility and optimality. Thus, the |Opt. gap| for NCW is zero. We follow the settings in [19] as the initial guess for IPOPT and SCIP. Based on the test results, the proposed HCM converges to the optimum of the distributed OEF problem with feasibility guarantee. Although references [3], [4], [7] derive promising results for PWL-, SOC- and CCP-based distributed OEF problems, they only consider a small number of agents ($\leq 4$). For this case that contains *more agents*, all their methods fail to converge within 1,000 seconds (considered infeasible), and the Obj. and |Opt. gap| are not applicable (N/A). We infer that the number of agents and binary variables may affect the convergence of PWL-, SOC-, and CCP-based methods.

TABLE I
COMPARISON BETWEEN DIFFERENT METHODS

| Method | PWL [3] | SOC [4] | CCP [7] | HCM | NCW |
|---|---|---|---|---|---|
| Time (s) | >1000 | >1000 | >1000 | 36.2 | 2.3 |
| Feasibility | ✗ | ✗ | ✗ | ✓ | ✓ |
| Obj. ($) | N/A | N/A | N/A | 3.714E+04 | 3.715E+04 |
| \|Opt. gap\| | N/A | N/A | N/A | 2.7E-04 | 0 |

Then, we test the scalability of the proposed HCM on the medium-scale IEGS-118-20, and Table II presents the results. Since PWL-, SOC-, and CCP-based methods do not converge for IEGS-6-7, we do not compare these methods with the proposed HCM anymore. This test indicates that the proposed HCM is competent to address distributed OEF problems in the IEGSs with a large number of agents, validating its scalability.

TABLE II
SCALABILITY OF THE PROPOSED HCM

| Method | Obj. ($) | Opt. gap | Feasibility |
|---|---|---|---|
| HCM | 4.207E+05 | 2.4E-04 | ✓ |
| NCW | 4.206E+05 | 0 | ✓ |

## V. CONCLUSION

This letter proposes a hidden convexity-based method to address multi-agent distributed OEF problems in transmission-level IEGSs. By exploring the hidden convexity, the nonconvex OEF subproblems are equivalently reformulated as SDPs. An ADMM is adopted to solve the hidden convexity-based distributed OEF problem. Test results indicate that compared with other methods (for handling Weymouth equations), the proposed hidden convexity-based method i) enhances the convergence of ADMM, and ii) exhibits its effectiveness in handling distributed OEF problems with *a large number of agents*. However, the distributed method may not completely guarantee the privacy of individual agents due to data interceptions, and we do not consider the multi-time scale problem in IEGSs. Future work includes i) combining distributed methods with encryption techniques to prevent potential interceptions, ii) addressing distributed multi-time scale OEF problems, and iii) extending the hidden convexity-based distributed method to multi-period OEF problems considering line packs and energy storage facilities.

## ACKNOWLEDGEMENT

This work adopts Mosek [20], Gurobi [21], IPOPT [22], and SCIP [23] as the general solvers for solving SDPs, mixed-integer convex programs, continuous convex and nonconvex programs, and mixed-integer nonconvex programs, respectively. The authors would like to thank all the researchers for developing these solvers. Meanwhile, the authors also greatly appreciate the editor and reviewers for helping improve the quality of this work.